\documentclass[a4paper,11pt]{amsart}

\usepackage[english]{my-shortcuts}
\usepackage{srcltx,algorithm,algorithmic}
\usepackage{enumitem}

\usepackage{tikz}
\begin{document}

\title[Von Neumann's inequality for tensors]
{Von Neumann's inequality for tensors}
\author{St\'ephane Chr\'etien{\tiny and} Tianwen Wei}

\address{Laboratoire de Math\'ematiques, UMR 6623\\ 
Universit\'e de Franche-Comt\'e, 16 route de Gray\\
25030 Besancon, France} 
\email{stephane.chretien@univ-fcomte.fr}

\address{Laboratoire de Math\'ematiques, UMR 6623\\ 
Universit\'e de Franche-Comt\'e, 16 route de Gray\\
25030 Besancon, France} 
\email{tianwen.wei@univ-fcomte.fr}

\begin{abstract}
For two matrices in $\mathbb R^{n_1\times n_2}$, the von Neumann inequality says that their scalar product is less than or equal to the scalar product of their singular spectrum. 
In this short note, we extend this result to real tensors and provide a complete study of the equality case.  
\end{abstract}

\maketitle

\section{Introduction}
The goal of this paper is to generalize von Neumann's inequality from matrices to tensors. Consider two matrices 
$X$ and $Y$ in $\mathbb R^{n_1\times n_2}$. Denote their singular spectrum, i.e. the vector of their 
singular values, by $\sigma(X)$ (resp. $\sigma(Y)$). The classical matrix von Neumann's inequality \cite{VonNeumann:TUR37} 
says that  
\bean 
\langle X,Y\rangle & \le & \langle \sigma(X),\sigma(Y)\rangle,
\eean 
and equality is achieved if and only if $X$ and $Y$ have the same singular subspaces. Von Neumann's 
inequality, and the characterization of the equality case in this inequality, are important in many 
aspects of mathematics. 

For tensors, the task of generalizing Von Neumann's inequality is rendered harder because of the necessity to 
appropriately define the singular values and the Singular Value Decomposition(SVD). In this paper, we will 
use the SVD defined in \cite{DeLathauwer:SIAMMatAnal00}, which is based on the Tucker decomposition. 

Our main result is given in Theorem \ref{VN3} below and gives a characterization of the equality case. We 
expect this result to be useful for the description of the subdifferential of some tensor fonctions as the
matrix counterpart has proved for matrix functions \cite{Lewis:JCA95}. Such functions occur naturally 
in computational statistics, machine learning and numerical analysis \cite{Gandy:InvProb11,Kressner:BIT13} 
due to the recent interest of sparsity promoting norms as a convex surrogate to rank penalization.    

\section{Main facts about tensors}
Let $D$ and $n_1,\ldots,n_D$ be positive integers. Let $\Xc \in \mathbb R^{n_1\times \cdots \times n_D}$ denote a 
$D$-dimensional array of real numbers. We will also denote such arrays as tensors. 
\subsection{Basic notations and operations}
A subtensor  of $\Xc$ is a tensor obtained by fixing some of its coordinates. As an example, fixing
one coordinate $i_d=k$ in $\Xc$ for some $k\in\{1,\ldots, n_d\}$ yields a tensor in
$\Rb^{n_1\times\cdots\times n_{d-1} \times n_{d+1}\times \cdots \times n_D}$.
 In the sequel, we will denote this subtensor of  $\Xc$ by $\Xc_{i_d=k}$.
 
The fibers of a tensor are subtensors that have only one mode, i.e. obtained by fixing every coordinate except one. The mode-$d$ fibers are the vectors
\bean 
\left(\Xc_{i_1,\ldots,i_{d-1},i_d,i_{d+1},\ldots,i_D}\right)_{i_d=1,\ldots,n_d}.
\eean 
They extend the notion of columns and rows from the matrix to the tensor framework. 
For a matrix, the mode-1 fibers are the columns and the mode-2 fibers are the rows. 

The mode-$d$
matricization $\Xc_{(d)}$ of $\Xc$ is obtained by forming the matrix whose columns are the mode-$d$ fibers 
of the tensor, arranged in the lexicographic ordering \cite{Kolda2009}. 
Clearly, the $k$th column of $\Xc_{(d)}$ consists of the entries of $\Xc_{i_d=k}$.

The mode-$d$ multiplication of a tensor $\Xc\in\Rb^{n_1\times\cdots \times n_D}$ by a matrix $U \in \Rb^{n_d'\times n_d}$, 
denoted by $\Xc\times_d U$, gives a tensor in $\Rb^{n_1\times\cdots \times n_d'\times\cdots \times n_D}$. It is  defined as 
\bean 
(\Xc\times_d U)_{i_1,\ldots,i_{d-1},i_d^\prime,i_{d+1},\ldots,i_D} 
& = & \sum_{i_d=1}^{n_d} \Xc_{i_1,\ldots,i_{d-1},i_d,i_{d+1},\ldots,i_D} U_{i_d',i_d}.
\eean 


\subsection{Higher Order Singular Value Decomposition (HOSVD)}
The Tucker decomposition of a tensor is a very useful decomposition, which can be chosen so that 
with appropriate 
orthogonal transformations, one can reveal a tensor $\Sc$ hidden inside $\Xc$ with interesting rank and orthogonality properties.
More precisely, we have 
\bea
\label{Tucker}
\Xc & = & \Sc(\Xc) \times_1 U^{(1)} \times_2 U^{(2)} \cdots \times_D U^{(D)}, 
\eea 
where each $U^{(d)}\in\Rb^{n_d\times n_d}$ is orthogonal and $\Sc(\Xc)$ is a tensor of the same size as $\Xc$ defined as follows. 
Moreover, subtensors $\Sc(\Xc)_{i_d=k}$ for $k=1,\ldots, n_d$ are all orthogonal to each other for each $d=1,\ldots,D$. 

\subsubsection{Relationship with matricization}
A tensor can be matricized along each of its modes. Let $\otimes$ denote the standard Kronecker product 
for matrices. Then, the mode-$d$ matricization of a tensor $\Xc \in 
\mathbb R^{n_1 \times \cdots \times n_D}$ is given by 
\bea
\label{matriz}  \label{1021b}
\Xc_{(d)} & = & U^{(d)} \cdot \Sc(\Xc)_{(d)}\cdot \left(U^{(d+1)} \otimes \cdots \otimes U^{(D)} \otimes U^{(1)} \otimes \cdots \otimes U^{(d-1)}\right)^t.  
\eea 
Take the (usual) SVD of the matrix $\Xc_{(d)}$
\bean 
\Xc_{(d)} & = & U^{(d)} \Sigma^{(d)} {V^{(d)}}^t
\eean  
and based on (\ref{matriz}), we can set  
\bean 
\Sc(\Xc)_{(d)} & = & \Sigma^{(d)} {V^{(d)}}^t\left(U^{(d+1)} \otimes \cdots \otimes U^{(D)} \otimes U^{(1)} 
\otimes \cdots \otimes U^{(d-1)}\right),
\eean 
where $\Sc(\Xc)_{(d)}$ is the mode-$d$ matricization of $\Sc(\Xc)$. One proceeds similarly for all $d=1,\ldots,D$
and one recovers the orthogonal matrices $U^{(1)},\ldots,U^{(D)}$ which allow us to decompose $\Xc$ as 
in (\ref{Tucker}). 
\subsubsection{The spectrum}
The mode-$d$ spectrum is defined as the vector of singular values of $\Xc_{(d)}$ and
we will denote it by $\sigma^{(d)}(\Xc)$. 
Notice that this construction implies that $\Sc(\Xc)$ has orthonormal fibers for every modes. With a slight abuse of notation, we will denote by $\sigma$ 
the mapping which to each tensor $\Xc$ assigns the vector $1/\sqrt{D}\: (\sigma^{(1)},\ldots,\sigma^{(D)})$ of all mode-$d$ singular spectra.

\section{Main Result}
\subsection{The main theorem}
The main result of this paper is the following theorem. 
\begin{theo}\label{VN3}
Let $\Xc,\Yc\in \Rb^{n_1\times \cdots \times n_D}$ be tensors.  
Then for all $d=1,\ldots,D$, we have 
\bea\label{vonneumann1}
\langle \Xc, \Yc\rangle \leq \langle \sigma^{(d)}(\Xc), \sigma^{(d)}(\Yc) \rangle.
\eea
The equality in (\ref{vonneumann1}) holds simultaneously for all $d=1,\ldots,D$ if and only there 
exist orthogonal matrices $W^{(d)}\in\Rb^{n_d\times n_d}$ for $d=1,\ldots, D$ 
and tensors $\Dc(\Xc),\Dc(\Yc)\in\Rb^{n_1\times\cdots\times n_D}$ such that
\bean
\Xc&=& \Dc(\Xc) \times_1 W^{(1)} \cdots \times_{D} W^{(D)}, \\
\Yc&=& \Dc(\Yc) \times_1 W^{(1)} \cdots \times_{D} W^{(D)},
\eean
where $\Dc(\Xc)$ and $\Dc(\Yc)$ satisfy the following properties:
 \begin{itemize}
 \item $\Dc(\Xc)$ and $\Dc(\Yc)$ are block-wise diagonal with the same number and size of blocks.
 
\item   Let $L$ be the number of blocks and $\{\Dc_l(\Xc)\}_{l=1,\ldots,L}$ (resp. $\{\Dc_l(\Yc)\}_{i=l,\ldots,L}$)
be the blocks on the diagonal of $\Dc(\Xc)$ (resp. $\Dc(\Yc)$).
  Then for each $l=1,\ldots, L,$ the two blocks $\Dc_l(\Xc)$ and $\Dc_l(\Yc)$ are proportional.
 \end{itemize}
 \end{theo}

\begin{figure}
\begin{center}\begin{tikzpicture}[scale=0.70]

\draw[thick] (0,0) -- (0,8) -- (8,8)--(8,0)--(0,0);
\draw[thick] (8,0) -- (12,2) -- (12,10) -- (4,10) -- (0,8);
\draw[thick] (8,8) -- (12,10);
\draw[dashed] (4,2) -- (12,2);
\draw[dashed] (4,2) -- (4,10);
\draw[dashed] (4,2) -- (0,0);

\draw (0,7) -- (2,7) -- (2,8) -- (3, 8.5);
\draw (1, 8.5) -- (3, 8.5) -- (3, 7.5) -- (2, 7);


\draw (3, 7.5) -- (3, 4.5) -- (4.5,4.5) -- (4.5, 7.5) -- (3, 7.5);
\draw (3, 7.5) -- (4.5, 8.25) -- (6, 8.25) -- (6, 5.25) -- (4.5, 4.5);
\draw (4.5, 7.5) -- (6, 8.25);


\draw (6, 5.25) -- (6, 4.25) -- (7, 4.25) -- (7, 5.25) -- (6, 5.25);
\draw (7.5, 4.5) -- (7.5, 5.5) -- (6.5, 5.5);

\draw (6, 5.25) -- (6.5, 5.5);
\draw (7, 5.25) -- (7.5, 5.5);
\draw (7, 4.25) -- (7.5, 4.5);


\draw (7.5, 4.5) -- (7.5, 4) -- (10, 4) -- (10, 4.5) -- (7.5, 4.5);
\draw (7.5, 4.5) -- (7.75, 4.625) -- (10.25, 4.625)-- (10.25, 4.125)-- (10, 4);
\draw (10.25, 4.625) -- (10, 4.5);

\draw (10.25, 4.125) -- (10.25, 1.625) -- (11.25, 1.625) -- (11.25, 4.125) -- (10.25, 4.125);
\draw (10.25, 4.125) -- (11, 4.5) -- (12, 4.5) -- (11.25, 4.125);
\end{tikzpicture}
\end{center}
\caption{A block-wise diagonal tensor.}
\end{figure}
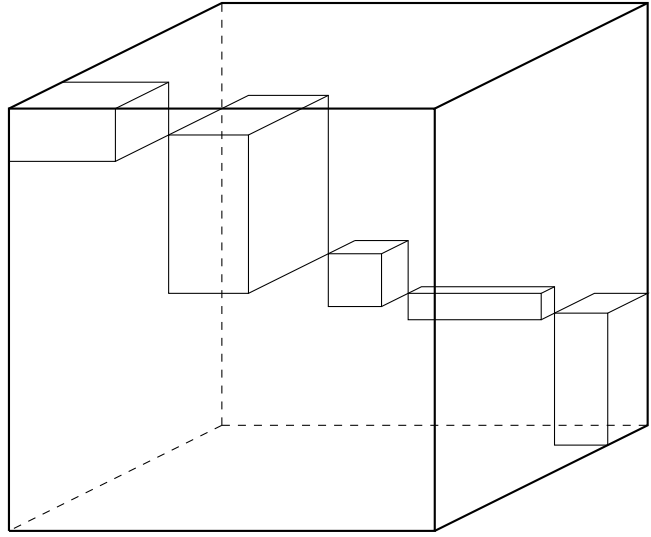

\subsection{Proof of the main theorem} 
In this section, we prove Theorem \ref{VN3}.
If $\Xc$ or $\Yc$ is a zero tensor, then the result is trivial. In the sequel, we assume that both
 $\Xc$ and $\Yc$ are non-zero tensors.
 
\subsubsection{The "if" part}
 The "if" part of the result is straightforward. Notice that 
$\langle \Xc, \Yc\rangle =\langle \Dc(\Xc), \Dc(\Yc)\rangle $ and
the singular vectors of $\Xc$ 
(resp. $\Yc$) are equal to those of $\Dc(\Xc)$ (resp. $\Dc(\Yc)$). Therefore, it remains to prove that
\bea\label{1103e}
\langle \Dc(\Xc), \Dc(\Yc)\rangle  = \langle \sigma^{(d)}( \Dc(\Xc)), \sigma^{(d)}(\Dc(\Yc)) \rangle,
\quad d=1,\ldots, D.
\eea
The conditions that $\Dc(\Xc)$ and $\Dc(\Yc)$ are block-wise diagonal and that $\Dc_i(\Xc)$ and $\Dc_i(\Yc)$ are proportional implies that each row of $\Dc_{(d)}(\Xc)$ and that of $\Dc_{(d)}(\Yc)$ are parallel. It then follows that $\Dc_{(d)}(\Xc)$ and $\Dc_{(d)}(\Yc)$ have the same left and right singular vectors. Then, applying the matrix
von Neumann's result immediately gives (\ref{1103e}).
 
\subsubsection{The "only if" part: first step}
Assume that
\bean
\langle \Xc, \Yc\rangle = \langle \sigma^{(d)}(\Xc), \sigma^{(d)}(\Yc) \rangle,\quad d=1,\ldots,D.
\eean
By the classical results of matrix von Neumann's inequality, we know that the equality holds if and only if
there exist orthogonal matrices $U^{(d)}$ and $V^{(d)}$ such that
\bea \label{1021a1}
\Xc_{(d)} & = & U^{(d)} {\rm Diag}\left(\sigma^{(d)}(\Xc)\right) {V^{(d)}}^t \quad \textrm{and}\\
\Yc_{(d)} & = & U^{(d)} {\rm Diag}\left(\sigma^{(d)}(\Yc)\right) {V^{(d)}}^t \nonumber
\eea
for all $d=1,\ldots, D$. From this remark, we obtain the following HOSVD of $\Xc$ and $\Yc$:
\bean 
\Xc & = & \Sc(\Xc) \times_1 U^{(1)} \cdots \times_{D} U^{(D)}, \\
\nonumber \\
\Yc & = & \Sc(\Yc) \times_1 U^{(1)} \cdots \times_{D} U^{(D)}. 
\eean

\subsubsection{The "only if" part: second step} 
We now show that subtensors $\Sc(\Xc)_{i_d=k}$ and $\Sc(\Yc)_{i_d=k}$ must be parallel
for all $k=1,\ldots, n_d$ and $d=1,\ldots,D$. 

Comparing (\ref{1021a1}) with (\ref{matriz}), we deduce that
\bea
\Sc_{(d)}(\Xc) & = &
\begin{pmatrix}
\sigma^{(d)}_1(\Xc) \cdot p_1^t \\
\vdots \\
\sigma^{(d)}_{n_d}(\Xc) \cdot p_{n_d}^t \\
\end{pmatrix},    \label{1021d1} 
\eea
where $p^t_i$ denotes the $i$th row of matrix ${V^{(d)}}^t\left(U^{(d+1)} \otimes \cdots \otimes U^{(D)} \otimes U^{(1)}  \otimes \cdots \otimes U^{(d-1)}\right)$.
Similarly, we have
\bea
\Sc_{(d)}(\Yc) =\begin{pmatrix}
\sigma^{(d)}_1(\Yc) \cdot p_1^t \\
\vdots \\
\sigma^{(d)}_{n_d}(\Yc) \cdot p_{n_d}^t \\
\end{pmatrix}.   \label{1021d2}
\eea
Comparing now (\ref{1021d1}) and (\ref{1021d2}) reveals that the $i${th} row of $\Sc_{(d)}(\Xc)$ and the $i${th} row of $\Sc_{(d)}(\Yc)$ must be proportional, for all 
$i=1,\ldots, n_d$. Formally, this means
\bea\label{1102a}
\sigma^{(d)}_{i_{d}}(\Yc) \cdot \Sc(\Xc)_{i_1 \cdots i_d  \cdots i_D}  
= \sigma^{(d)}_{i_{d}}(\Xc) \cdot  \Sc(\Yc)_{i_1 \cdots i_{d} \cdots i_D}.
\eea
for all possible values of $i_1,\ldots, i_D$.

\subsubsection{The "only if" part: third step} For $d=1,\ldots,D$, let $r^{(d)}_x$ (resp. $r^{(d)}_y$) be the rank of 
$\Sc_{(d)}(\Xc)$ (resp. $\Sc_{(d)}(\Yc)$). Let $(i_1,\ldots,i_D)$ be such that 
\bean 
(i_1,\ldots,i_D) & \not \le & (r^{(1)}_y,\ldots,r^{(D)}_y) \\
(i_1,\ldots,i_D) & \not \ge & (r^{(1)}_y,\ldots,r^{(D)}_y).
\eean 
Then, $\Sc(\Yc)_{i_1,\ldots,i_D}=0$ but there exists $d\in \{1,\ldots,D\}$ such that $\sigma^{(d)}_{i_d}(\Yc)>0$. Using 
(\ref{1102a}), we obtain that $\Sc(\Xc)_{i_1,\ldots,i_D}=0$. Thus, if $r_x^{(d)}>r_y^{(d)}$ for some 
$d=1,\ldots,D$, then there exists some  
\bean 
(i_1,\ldots,i_D) & > & (r^{(1)}_y,\ldots,r^{(D)}_y)
\eean 
such that $\Sc(\Xc)_{i_1,\ldots,i_D}\neq 0$ and thus, $r_x^{(d)}>r_y^{(d)}$ for all $d=1,\ldots,D$. By symmetry, 
we deduce that either $r_x^{(d)}>r_y^{(d)}$ for all $d=1,\ldots,D$, or $r_x^{(d)}<r_y^{(d)}$ for all $d=1,\ldots,D$ or 
else $r_x^{(d)}=r_y^{(d)}$ for all $d=1,\ldots,D$.

\subsubsection{The "only if" part: fourth step} Assume that $(r_x^{(1)},\ldots,r_x^{(D)})\le (r_y^{(1)},\ldots,r_y^{(D)})$. 
The other case may be treated in the same way (with an overlap in the equality case) by interchanging the role of $\Xc$ and $\Yc$. 
For all $(i_1,\ldots,i_D)\le (r_y^{(1)},\ldots,r_y^{(D)})$, we have $\sigma^{(d)}_{i_d}(\Yc)>0$ for all $d=1,\ldots,D$. 
Thus, (\ref{1102a}) gives 
\bean
\Sc(\Xc)_{i_1 \cdots i_d  \cdots i_D}  
& = & \frac{\sigma^{(d)}_{i_{d}}(\Xc)}{\sigma^{(d)}_{i_{d}}(\Yc) } \cdot  \Sc(\Yc)_{i_1 \cdots i_{d} \cdots i_D}.
\eean
We deduce from this equation that for two indices $(i_1,\ldots,i_D)$ and $(i^\prime_1,\ldots,i^\prime_D)$, 
if 
\begin{itemize}
\item[(i)] there exists some $d$ in $\{1,\ldots,D\}$ such that
 $i_d=i^\prime_d$, 
 \item[(ii)] $\Sc(\Xc)_{i_1, \cdots, i_d,  \cdots, i_D}$ 
and $\Sc(\Xc)_{i_1^\prime, \cdots, i_d^\prime,  \cdots, i_D}$ are different from zero, 
\end{itemize}
then
\bean
(\Sc(\Xc)_{i_1 \cdots i_d  \cdots i_D},\Sc(\Xc)_{i_1^\prime, \cdots, i_d^\prime,  \cdots, i_D})
& = & \rho \cdot
(\Sc(\Yc)_{i_1 \cdots i_d  \cdots i_D},\Sc(\Yc)_{i_1^\prime, \cdots, i_d^\prime,  \cdots, i_D}),
\eean  
where  
\bean 
\rho=\frac{\sigma^{(1)}_{i_{1}}(\Xc)}{\sigma^{(1)}_{i_{1}}(\Yc) } & = \cdots = & \frac{\sigma^{(D)}_{i_{D}}(\Xc)}{\sigma^{(D)}_{i_{D}}(\Yc) } >0.
\eean  

\subsubsection{The "only if" part: fifth step} Let $\rho_1>\cdots>\rho_L$ denote the possible values of the ratio $\sigma^{(d)}_{i_{d}}(\Xc)/\sigma^{(d)}_{i_{d}}(\Yc)$, 
for all $(i_1,\ldots,i_D)\le (r_y^{(1)},\ldots,r_y^{(D)})$. Let $I_{d,l}$, $d=1,\ldots,D$, $l=1,\ldots,L$ 
denote the possibly empty set of indices in $\{1,\ldots,r_x^{(d)}\}$ such that 
\bean 
\frac{\sigma^{(d)}_{i_{d}}(\Xc)}{\sigma^{(d)}_{i_{d}}(\Yc) } & = & \rho_l
\eean 
and let $m_{d,l}$ denote the cardinality of $I_{d,l}$. Then, for each $d=1,\ldots,D$, we can find a permutation $\pi_d$ 
on $\{1,\ldots,n_d\}$ such that $\pi_d(I_{d,1})=\{1,\ldots,m_{d,1}\}$, $\pi_d(I_{d,2})=\{m_{d,1}+1,\ldots,m_{d,1}+m_{d,2}\}$, 
and so on and so forth. 

Thus, for each mode $d=1,\ldots,D$, there exists a permutation matrix $\Pi_d$ such that the matrices 
\bean
\Dc(\Xc) &=& \Sc(\Xc) \times_1 \Pi_1 \cdots \times_D \Pi_D, 
\eean 
and 
\bean 
\Dc(\Yc) &=& \Sc(\Yc) \times_1 \Pi_1 \cdots \times_D \Pi_D.
\eean
contain $L$ blocks and each block in $\Dc(\Xc)$ is proportionnal to 
the corresponding block in $\Dc(\Yc)$. Moreover, any entry $\Dc(\Xc)_{i_1,\ldots,i_D}$
with $(i_1,\ldots,i_D) \le (r_x^{(1)},\ldots,r_x^{(D)})$ and lying outside the union of these $L$ blocks 
is null since if it were not, by (\ref{1102a}) combined with $(r_x^{(1)},\ldots,r_x^{(D)})\le (r_y^{(1)},\ldots,r_y^{(D)})$, 
it would be proportionnal to a nonzero component $\Dc(\Yc)_{i_1,\ldots,i_D}$ with two different ratios, 
thus a contradiction. 

Finally, setting $W^{(d)}=\Pi_d \: U^{(d)}$ for $d=1,\ldots,D$ achieves the proof.

\bibliographystyle{amsplain}
\bibliography{database}

\end{document}